\documentclass[11pt,a4paper]{article}

\newcommand{\marginnote}[1]%
       {\mbox{}\marginpar{\raggedleft\hspace{0pt}\itshape #1}}%
\newcommand{\ctitlefont}%
   {\fontencoding{OT1}%
    \fontfamily{cmr}%
    \fontseries{bc}%
    \fontshape{n}%
    \fontsize{35.83}{54}%
    \selectfont}%
\newcommand {\beq  } {\begin{equation}     }%
\newcommand {\eeq  } {\end{equation}       }%
\newcommand{\be}[1]{\begin{equation}\label{eq:#1}}
\newcommand{\ee}{\end{equation}}
\newcommand\beds{\begin{description}}
\newcommand\eds{\end{description}}
\newcommand\benn{\begin{displaymath}}
\newcommand\eenn{\end{displaymath}}
\newcommand\bea[1]{\begin{eqnarray}\label{eq:#1}}
\newcommand\eea{\end{eqnarray}}

\newcommand\bfg[1]{\begin{figure}[h]\vspace{#1cm}}
\newcommand\efg{\end{figure}}

\newcommand\btb{\begin{table}}
\newcommand\etb{\end{table}}

\input psfig.sty
\newcommand{\Titolo}{\vskip 0.0cm \fontsize{24pt}{24pt} \selectfont}
\newcommand{\Abstract} {
     \vskip 0.3cm \fontsize{12pt}{12pt} \selectfont
     Abstract
     \vskip 0.3cm
     \fontsize{10pt}{10pt}
     \selectfont
     }
\newcommand{\Report}{\vskip 0.4cm \fontsize{11pt}{11pt} \selectfont}
\newcommand{\myline}{ \noindent \rule{13cm}{0.75pt} }
\newcommand{\Autore}{\vskip 0.0cm \fontsize{12pt}{12pt} \selectfont }
\newcommand{\Indirizzo}{\vskip 0.0cm \fontsize{10pt}{10pt} \selectfont }

\begin{document}

\Titolo \centerline{Simulations of Gaussian Processes}
\centerline{and Neuronal Modeling}

\myline \Autore {\noindent Elvira Di Nardo} \Indirizzo {
Dipartimento di Matematica, Universit\`a della Basilicata,
Contrada Macchia Ro-} \Indirizzo {mana, Potenza. E-mail:
dinardo@unibas.it}

\Autore {\noindent Amelia G. Nobile} \Indirizzo { Dipartimento di
Matematica e Informatica, Universit\`a di Salerno, Via S.
Allende,} \Indirizzo { Salerno. E-mail: nobile@unisa.it }

\Autore {\noindent         Enrica Pirozzi and Luigi M. Ricciardi}
\Indirizzo { Dipartimento di Matematica e Applicazioni,
Universit\`a di Napoli Federico II, Via } \Indirizzo {Cintia,
Napoli. E-mails: \{enrica.pirozzi, luigi.ricciardi\}@unina.it}

\myline

\Abstract
%% ======================================================================= %%
%% ABSTRACT ABSTRACT ABSTRACT ABSTRACT ABSTRACT ABSTRACT ABSTRACT ABSTRACT %%
%% ======================================================================= %%
{\it 
The research work outlined in the present note highlights the essential role played by the simulation procedures implemented by us on CINECA supercomputers to complement the mathematical investigations carried within our group over the past several years. The ultimate target of our research is the understanding of certain crucial features of the information processing and transmission by single neurons embedded in complex networks. More specifically, here we provide a bird's eye look of some analytical, numerical and simulation results on the asymptotic behavior of first passage time densities for Gaussian processes, both of a Markov and of a non-Markov type. Several figures indicate significant similarities or diversities between computational and simulated results.
}

 \Report
%% ======================================================================= %%
%% REPORT REPORT REPORT REPORT REPORT REPORT REPORT REPORT REPORT REPORT   %%
%% ======================================================================= %%
%\baselineskip 1cm

\section{Introduction to FPT models for the neuronal firing}
%INSERIMENTO DELLA RICERCA NEL QUADRO DISCIPLINARE:

The research work outlined here takes place within
the framework of applied probability. Our aim is to describe the dynamics of the neuronal 
firing by modeling it via a
stochastic process representing the change in the neuron membrane
potential between each pair of consecutive spikes (cf., for instance,
\cite{Ricciardi77}). In our approach, the threshold voltage is
viewed as a deterministic function, and the instant when the
membrane potential reaches it (i.e. when a spike occurs) as a
first passage time (FPT) random variable. We shall focus our
attention on neuronal models rooted on Gaussian processes,
partially motivated by the generally accepted hypothesis that in
numerous instances the neuronal firing is caused by the
superposition of a very large number of synaptic input pulses
which is suggestive of the generation of Gaussian distributions
 by virtue of some sort of central limit theorems.
\par
Let us first consider a one-dimensional non-singular Gaussian
stochastic process $\left\{X(t), t\geq  t_0  \right\}$ and a boundary
 $S(t)\in C^1[t_0,+\infty)$.
%\par{\bf Conditioned FPT model}
 We assume
$P\{X(t_0)=x_0\}=1$, with $x_0 < S(t_0),$ i.e. we focus our attention on the subset of sample
paths of $X(t)$ that originate at a preassigned state $x_0$ at the
initial time $t_0$. Then,
$$
T_{x_0} = \inf_{t\geq t_0} \bigl\{t : X(t) > S(t)\bigr\},\qquad
x_0<S(t_0)
$$
is the FPT of $X(t)$ through $S(t)$, and
\begin{equation}
g(t|x_0,t_0) = \frac{\partial}{\partial t} P(T_{x_0} < t)
\label{gcond}
\end{equation}
is its probability density function (pdf).
\par
Henceforth, the FPT pdf  $g(t|x_0,t_0)$
will be identified with the firing pdf of a neuron whose membrane
potential is modeled by $X(t)$ and whose firing threshold is
$S(t)$.
\par
%{\bf Upcrossing FPT model} \par \noindent

 In order to  include
more physiologically significant features -- such as a finite
decay constant of the membrane potential, the presence of reversal potentials,
time-dependent firing thresholds -- and to refer to wider classes
of inputs as responsible for the observed sequences of output
signals released by the neuron, we also define the FPT upcrossing
model. This is viewed as an FPT problem to a threshold, or boundary,
$S(t)$ for the subset of sample paths of the one-dimensional
non-singular Gaussian process $\left\{X(t), t \geq t_0 \right\}$
originating at a
state $X_0.$ Such initial state, in turn, is viewed as a random variable with pdf
%
%\begin{equation}
$$
\gamma_{\varepsilon}(x_0,t_0) \equiv \left\{\begin{array}{ll}
\displaystyle{ {f(x_0)\over P\{X(t_0)<S(t_0)-\varepsilon\}}}\,, \,
& x_0 < S(t_0) - \varepsilon \\
\hfill \\
0,  &  \, x_0  \geq S(t_0) -  \varepsilon.
\end{array} \right.
 \label{(gamma)}
$$
%\end{equation}
%
Here, $\varepsilon > 0$ is a fixed real number and $f(x_0)$
denotes the Gaussian pdf of $X(t_0)$. Then,
$$
T_{X_0}^{(\varepsilon)} = \inf_{t\geq t_0} \{t : X(t) > S(t)\},
$$
is the $\varepsilon$-upcrossing FPT of $X(t)$ through $S(t)$ and
the related pdf is given by
%
%\begin{equation}
$$
g_u^{(\varepsilon)}(t|t_0) = \frac{\partial}{\partial t}
P(T_{X_0}^{(\varepsilon)} < t) =\int_{-\infty}^{S(t_0)-\varepsilon}
\!\!\!\! g(t|x_0,t_0)\,\gamma_{\varepsilon}(x_0,t_0)\; dx_0\qquad (t\geq
t_0), \label{(updens)}
$$
%\end{equation}
%
where
%\begin{equation}
$
g(t|x_0,t_0)
$
%\end{equation}
%
is defined in (\ref{gcond}).
Without loss of generality, we set $t_0=0$ and $x_0=0$ and for this case we write $g(t) \equiv g(t|0,0)$ and ${g}_u(t) \equiv g_u^{(\varepsilon)}(t|0)$, for fixed values of $\epsilon.$ 

\par
%\section{Breve rassegna of the preliminary results}
The selection of one of the various methods available  to compute the
firing pdf's $g(t)$ and $g_u(t)$ depends on
the assumptions made on $X(t)$. For diffusion
processes (cf. \cite{Buonocore87}, \cite{Dicrescenzo00}) and for
Gauss-Markov processes (cf. \cite{a.a.p.}) we have proved that the
firing pdf is solution of a second kind integral Volterra equation.
For generally
regular thresholds we have designed, and successfully implemented, a
fast and accurate numerical procedure for solving such integral equation, and compared our approximations with those stemming out of standard numerical methods.
Furthermore, in \cite{a.a.p.}, by adopting a symmetry-based approach, we have determined the exact firing pdf for
thresholds  of a suitable analytical form.  

Mathematical models based on non-Markov stochastic processes better describe the correlated firing activity, even though their analytical treatment is more complicated and only rare and
fragmentary results appear to be available in literature. 
\par
By using a variant of the method proposed in \cite{Ricciardi4}, for a zero-mean non-singular stationary Gaussian process differentiable in the mean square sense, a cumbersome series expansion for the conditioned
FPT pdf density and for the upcrossing FPT pdf
density has been obtained (\cite{vienna2002}, \cite{kluwer}). In both cases we have succeeded to obtain numerically a reliable
evaluation of the first term of the series expansion. By comparisons of these results with the simulated firing densities (\cite{iscs99},
\cite{vienna2002}), we
have been led to conclude that this first term is a good approximation of firing pdf's only
for small times.

The simulation procedures provide valuable alternative
investigation tools
 especially if they can
be implemented on parallel computers, (see \cite{DiNardo01}). We point out that our simulation
originates from Franklin's algorithm
\cite{Franklin}. We have implemented it in both vector and
parallel modalities after
suitably modifying it for our computational needs, for instance to obtain reliable approximations of upcrossing densities
(cf. \cite{DiNardo4}, \cite{iscs99}). Thus doing,
reliable histograms of FPT densities of stationary Gaussian
processes with rational spectral densities can be obtained in the
presence of various types of boundaries.
\par
 For a thorough description
of the simulation algorithm we refer the reader to
\cite{DiNardo3}, \cite{eur99}, \cite{iscs99}.

We wish to stress that our endeavors strive to improve the simulation techniques, particularly within the context of stationary Gaussian processes for which alternative simulation procedures are being implemented and tested. Besides the use of algorithms based on special properties of the spectral density, such as its being of a rational type, methods based on the spectral representation of the process and circulant embedding methods are presently under investigation. The aim is to design efficient algorithms to simulate Gaussian processes of a more general type.
It must also be pointed out that within our research project the role of the simulation procedure is threefold:

\begin{description}
\item{(i)  } to provide an investigation tool to explore the behavior of the firing density in a variety of different conditions;
\item{(ii) } to permit us to evaluate reliability and precision of the results obtained via numerical and analytic approximations;
\item{(iii)} to represent the only possible alternative whenever analytic and computational methods are failing.
\end{description}
\par
\begin{figure}[t]
\centerline{
\psfig{figure=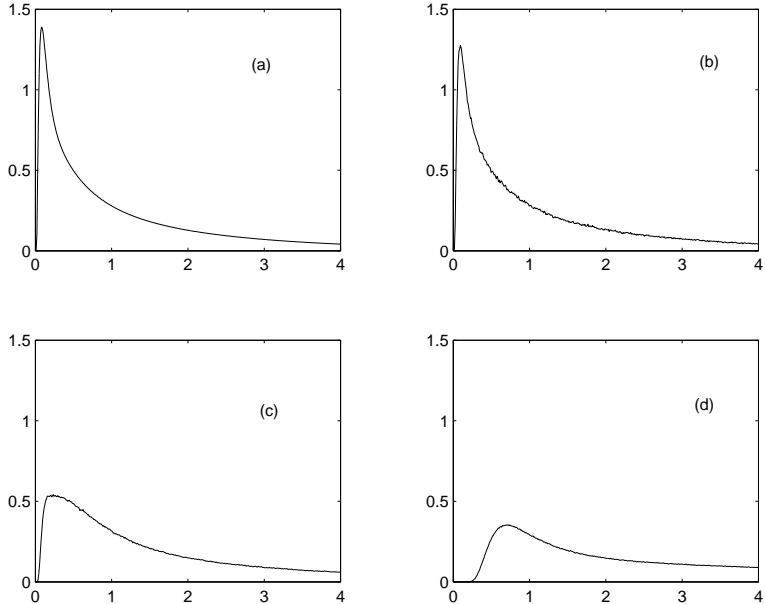,height='350pt}} 
\caption{ \small{
Plots refer to  FPT pdf $g(t)$ with $\beta=0.5$ for a zero-mean Gaussian process
            characterized by correlation function (\ref{(funcorr)}) in the presence of the boundary (\ref{soglia}) with $d=0.25$.
            In Figure~1(a), the function $g(t)$  has been plotted.
            The estimated FPT pdf $\tilde{g}(t)$ with $\alpha=10^{-10}$ is shown in Figure~1(b),
            with $\alpha=0.25$ in Figure~1(c) and with $\alpha=0.5$ in Figure~1(d).}} \label{fig:fig1}
\end{figure}

\section{Markov and non-Markov models: an analysis by simulations}
In order to analyze how the lack
of memory affects the shape of the FPT densities, in \cite{biocomp}  we have compared the behavior of such densities for Gauss-Markov processes versus Gauss non-Markov
processes. 
\par
Let us consider a zero-mean stationary Gaussian
process $X(t)$ with correlation function
\begin{equation}
\gamma(t)= e^{-\beta\, |t|}\, \cos (\alpha\, t), \qquad \alpha, \beta \in
{ I\!\! R}^{+} \label{(funcorr)}
\end{equation}
which is the simplest type of correlation of concrete interest for applications \cite{Yaglom87}.
Furthermore, let us assume that the boundary is of the following type:

\begin{equation}
S(t)=d\,e^{-\beta\, t}\, \Biggl\{ 1 - \frac{e^{2\, \beta\, t}-1}{2\, d^2}\, 
\ln \Biggl[ {1\over 4}+{1\over 4}\,
\sqrt{1+8 \exp\biggl(-{\displaystyle\frac{4\,d^2}{e^{2\,\beta\, t} -1}\biggr)}}
\;\;\Biggr]\Biggr\},
\label{soglia}
\end{equation}
with $d>0$.
Due to the assumed correlation
(\ref{(funcorr)}), $X(t)$ is  not mean square
differentiable (see \cite{biocomp} for the details). Thus, the afore-mentioned
series expansion (see \cite{kluwer}) does not hold. However,
specific assumptions on the parameter $\alpha$ help us to
characterize the shape of the FPT pdf.

We start assuming $\alpha=0,$ so that the correlation function
(\ref{(funcorr)}) factorizes as
$$
\gamma(t)= e^{-\beta\, \tau}\, e^{-\beta\,(t-\tau)} \qquad \beta
\in { I\!\! R}^{+}, \, 0<\tau<t.
$$
In such case $X(t)$ becomes  Gauss-Markov. 
As proved in \cite{a.a.p.}, for the Gauss-Markov process $X(t)$ with covariance (\ref{(funcorr)}), the FPT pdf $g(t)$ in the presence of boundaries of type (\ref{soglia}) can be evaluated in closed form. Alternatively, $g(t)$ can be numerically obtained by solving a non-singular Volterra integral equation (see \cite{a.a.p.}).
In Figure~1(a) such density is plotted for $\beta=0.5$ and
$d=0.25$.
\par
Setting $\alpha \neq 0$ in (\ref{(funcorr)}), the Gaussian process
$X(t)$ is no longer Markov and its spectral density is given by
\begin{equation}
\Gamma(\omega) = \frac{2\, \beta\,
(\omega^2+\alpha^2+\beta^2)}{\omega^4
+2\,\omega^2\,(\beta^2-\alpha^2)+(\beta^2+\alpha^2)^2}\,,
\label{(spectral)}
\end{equation}
thus being of a rational type. Since in (\ref{(spectral)}) the
degree of the numerator is less than the degree of the
denominator, it is possible to apply the simulation algorithm
described in \cite{eur01} in order to estimate the FPT pdf
$\tilde{g}(t)$ of the process.
\par
The simulation procedure has been implemented by a parallel
FORTRAN 90 code on a 128-processor IBM SP4 supercomputer, based on
MPI language for parallel processing, made available to us by CINECA. 
The number of simulated
sample paths has been set equal to $10^7$. The estimated FPT pdf's
$\tilde{g}(t)$ through the specified boundary are plotted in
Figures~1(b)$\div$1(d) for Gaussian processes with correlation
function (\ref{(funcorr)}) in which we have taken $\alpha=10^{-10}, 0.25, 0.5$,
respectively. Note that as $\alpha$ increases, the shape of the
FPT pdf $\tilde{g}(t)$ becomes progressively flatter, while its mode
increases. Furthermore, as Figures 1(a)-1(b) indicate, $\tilde{g}(t)$ is
very close to $g(t)$ for small values of $\alpha$.

%------------
% FIGURE 2
%------------

\begin{figure}[t]
\centerline{
\psfig{figure=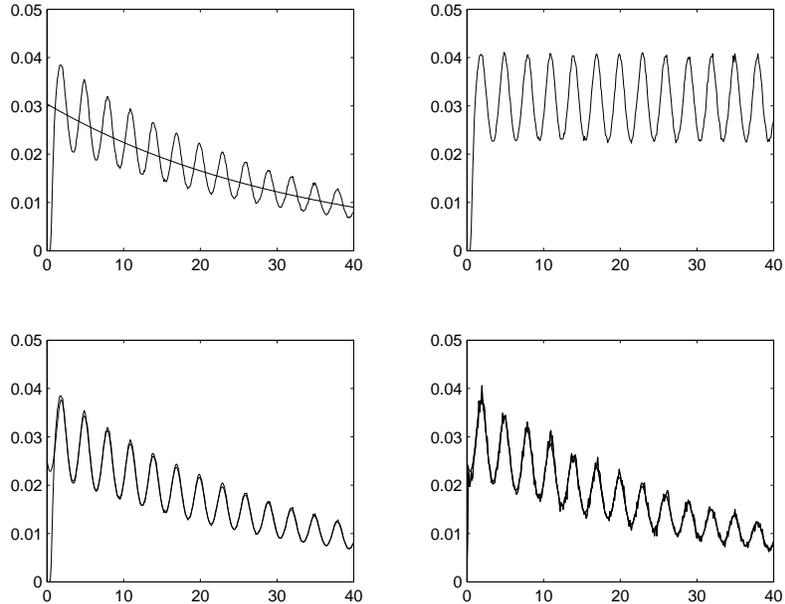,height='350pt}} 
\caption{ {\small
For a zero-mean stationary Gaussian process with
covariance (\ref{es:gencovariance}) and $\alpha=\beta=1$ in the presence
of the periodic boundary $S(t)=2+0.1\,\sin(2\,\pi\,t/3),$ the simulated FPT density ${\tilde g}(t)$ is compared with the 
exponential density ${\widehat\lambda}\, e^{-{\widehat\lambda}\, t}$, with   
${\widehat\lambda}=0.030386$ in (a). The function ${\tilde Z}(t)=\tilde{g}(t)\, 
\exp\bigl\{{{\widehat\lambda}\,t}\bigr\}$ is plotted in (b). 
The asymptotic
exponential approximation is plotted together with the simulated
FPT density $\tilde{g}(t)$ in (c).
The asymptotic
exponential approximation  is plotted together with the simulated
upcrossing FPT density $\tilde{g}_u(t)$ in (d).}} \label{fig2report}
\end{figure}

\section{Asymptotic results}

The asymptotic behavior of the FPT densities
for Gaussian processes as boundaries or time grow larger has been
studied in \cite{eur01}, \cite{vienna2002} and  \cite{kluwer}. Our 
analysis is a natural extension of some
investigations performed for the Ornstein-Uhlenbeck (OU) process
\cite{Nobile1} and successively extended to the class of
one-dimensional diffusion processes admitting steady state
densities in the presence of single asymptotically constant
boundaries or of single asymptotically periodic boundaries
(see \cite{Giorno} and \cite{Nobile}). There, computational as well as
analytical results have indicated that the conditioned FPT pdf is
susceptible of an excellent non-homogeneous exponential
approximation for large boundaries, even though these boundaries need not be very distant
from the initial state of the process. To this aim, we have
estimated such a density by generating the sample paths of the
Gaussian process through the parallel simulation algorithm
implemented on the super-computer CRAY T3E of CINECA in order to overcome
the outrageous complexity offered by the numerical evaluation of
the involved partial sums of the conditioned FPT pdf series
expansion. Specifically, we have considered the class of zero-mean stationary
Gaussian processes characterized by damped oscillatory covariances
(see \cite{eur01}):

\begin{eqnarray}
\gamma (t)= e^{-\beta\,|t|} \,\Bigl[\cos\bigl(\alpha\, t)+
\,\sin\bigl(\alpha |t|)\Bigr],
\label{es:gencovariance}
\end{eqnarray}
where $\alpha$ and $\beta$ are positive real numbers.

\par
 %After the analysis of the case of asymptotically
%constant thresholds (see, for example, Figure \ref{f:fig2}), 
The
results of our computations have shown that  for certain
periodic boundaries of the form
\begin{equation}
S(t)=S_{0} + A\,\sin \bigl(2\,\pi\,t/Q\bigr), \quad S_{0}, A, Q >0,
\label{boundary}
\end{equation}
not very distant from the initial value of the process, the FPT
pdf soon exhibits damped oscillations having the same period of
the boundary. Furthermore, starting
from rather small times, the estimated FPT densities $\tilde{g}(t)$
appears to be representable in the form
\begin{equation}
\tilde{g}(t) \simeq \tilde{Z}(t)\; e^{-\widehat{\lambda}\, t},
\label{stima}
\end{equation}
where $\widehat{\lambda}$ is a constant that can be estimated  by
the  least squares methods, and where $\tilde{Z}(t)$ is a periodic
boundary of  period $T$ (see, for example, Figure \ref{fig2report}(a)
and Figure \ref{fig2report}(b)).
The goodness of the exponential approximation increases as the
boundary is progressively moved farther apart from the starting
point of the process. The more the periodic boundary is far from the
starting point of the process, the more the exponential approximation
improves.
\par
In \cite{kluwer} we have shown by rigorous mathematical arguments that as boundary  
(\ref{boundary}) moves away from the initial state of the process, the FPT pdf approaches a non-homogeneous exponential density of the type
\begin{equation}
\tilde{g}(t) \sim \tilde{h}(t) \exp \left\{ -
\int_0^t \tilde{h}(\tau) d\tau \right\}, \label{eq:1}
\end{equation}
where 

\begin{eqnarray}
&&\tilde{h}(t)=
{\sqrt{\alpha^2+\beta^2}\over 2\,\pi} \exp \biggl\{ - {S^2(t)\over 2} \biggr\}\, 
\biggl[ \exp \biggl\{ - {[\dot{\rho}(t)]^2\over 2\,(\alpha^2+\beta^2)} 
\biggr\}\nonumber \\ 
&&\hspace*{4cm}-  \sqrt{{\pi\over 2\,(\alpha^2+\beta^2)}} \, \dot{\rho}(t) \, {\rm Erfc} 
\, \biggl({\dot{\rho}(t)\over\sqrt{2\,(\alpha^2+\beta^2)}}  \biggr) \biggr]
\label{es:R[Z(t)]}
\end{eqnarray}
and $\rho(t)=A \sin (2\pi t/Q).$ (See, for example, Figure \ref{fig2report}(c)).
In \cite{vienna2002} a similar result is proved for the
upcrossing FPT density. (See Figure \ref{fig2report}(d)).

%***************************
% Figure 3
%***************************
\begin{figure}[t]
\centerline{\psfig{figure=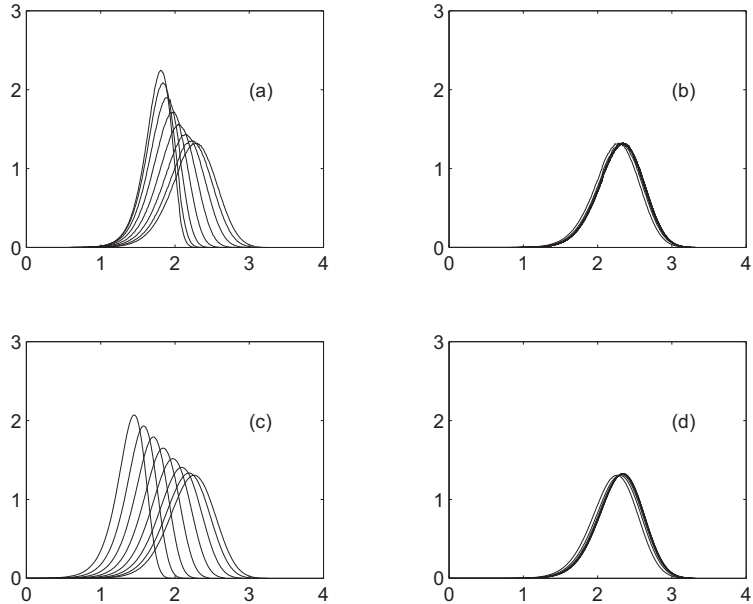,height='350pt}}
  \caption{ {\small  For different choices of $\vartheta$ in the interval $[0.008,1.024]$, with threshold $S(t) = -t^2/2
- t+5$, plot of the simulated
$\tilde{g}_u(t)$ is shown in (a) and plot of
$\tilde{g}_u(t)$ for the OU-model in (c). Same in (b) and (d)
for values of $\vartheta$  in $[2.048,200].$}} \label{f:fig5}
\end{figure}
%
%***************************
% Figure 4
%***************************

%
\begin{figure}[h]
\centerline{\psfig{figure=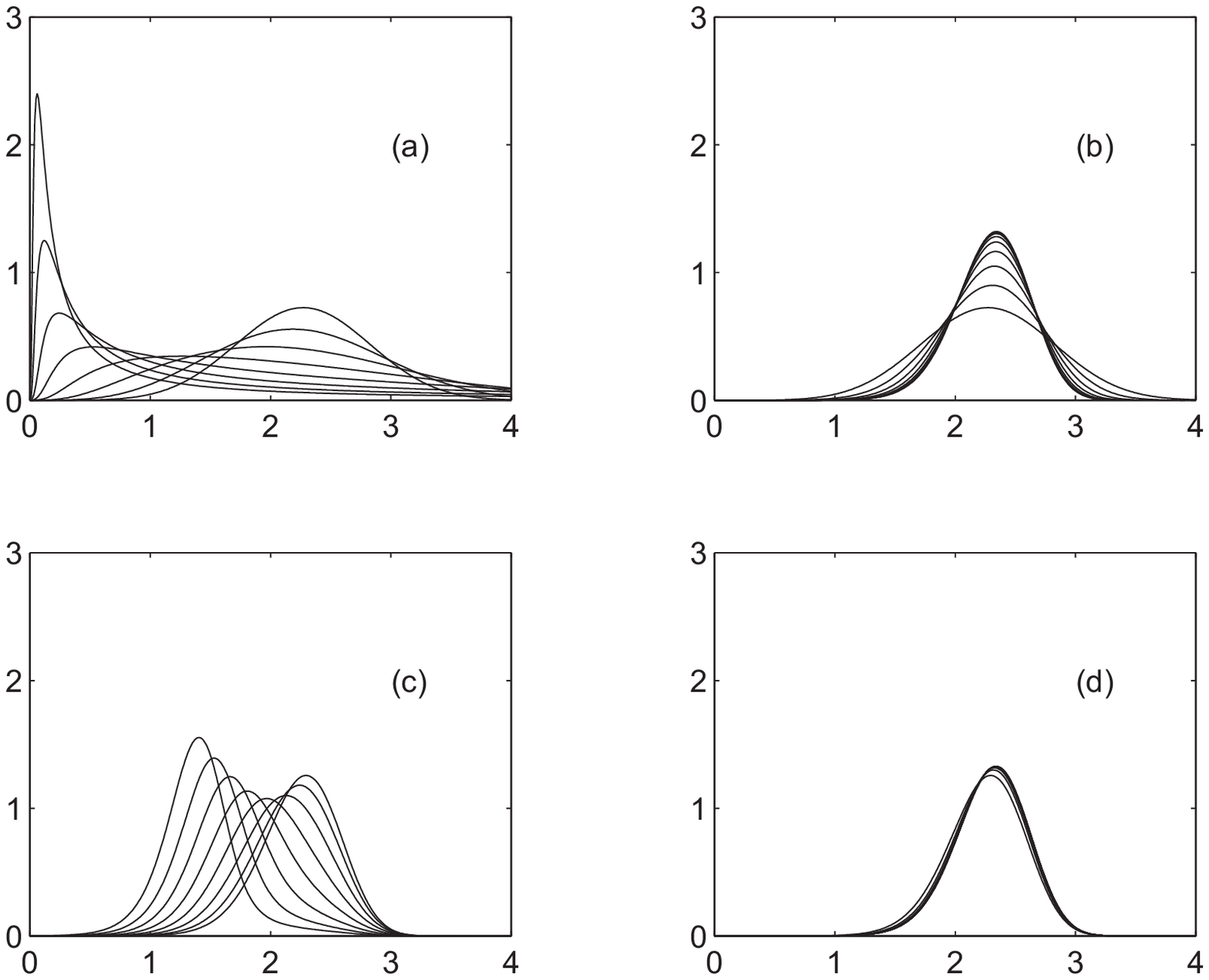,height='350pt}}
\caption{{ \small For different choices of $\vartheta$ in the interval $[0.008,1.024]$, with the same boundary as in Figure \ref{f:fig5}, plot of $\tilde{g}_u(t)$ for the Wiener model is shown in
(a) and plot of $q(t)$ for the
Kostyukov-model in (c). Same in (b) and (d) for values of $\vartheta$  in $[2.048,200].$}}
\label{f:fig6}
\end{figure}

 It should be stressed that the
analytic and the simulation results constitutes only a preliminary
step towards the construction of neuronal models based on
non-Markov processes. Nevertheless, the unveiling of properties of
the asymptotic behavior of  FPT may turn out to be useful also for
the description of neuronal activities at small times whenever
the intrinsic time scale of the microscopic events involved during
the neuron's evolution is much smaller than the macroscopic
observation time scale, or when the asymptotic regime is exhibited
also in the case of firing thresholds not too distant from the
resting potential, similarly to what was already pointed out by us
in connection with the OU neuronal model
\cite{Giorno}.

%%%%%%%%%%%

\section{An alternative approach}
Within the context of single neuron's activity modeling a completely
different, apparently not well known, approach was proposed by
Kostyukov {\it et al.} (\cite{Kost81}) in which a non-Markov
process of a Gaussian type is assumed to describe the time course
of the neural membrane potential.
The model due to Kostyukov (K-model) makes use of the notion of correlation
time. Namely, let
 $X(t)$ be a stationary Gaussian process with zero mean, unit variance and
correlation function $R(t)$. Then \cite{Strato},
$$
\vartheta=\int_{0}^{+\infty} |R(\tau)|\; d \tau < +\infty
$$
is defined as the correlation time of the process $X(t).$  Under
some assumptions on the threshold and by using a sort of diffusion
approximation, Kostyukov works out a numerical
evaluation $q(t)$ to the upcrossing FPT pdf. This approximation is
obtained as solution of an integral equation that  can be solved
by routine methods. The relevant feature of this approach is that
the unique parameter $\vartheta$ characterizes the considered
class of stationary standard Gaussian processes.

In \cite{Dicrescenzo00} and in \cite{DiNardo4}  we analyzed this
method pinpointing similarities and differences with respect to our
models. Recently \cite{eur03}, we have again made use of the K-model 
and compared the obtained results with those worked out by  numerically
solving the integral equation holding for Gauss-Markov
processes in the case of the OU and Wiener models, see \cite{DiNardo4}
for details, and with the results obtained via the simulations of Gaussian
processes. A variety of thresholds and of $\vartheta$
values has been considered. For some values of $\vartheta$
between 0.008 and 200, our results have been plotted in Figures \ref{f:fig5} and
\ref{f:fig6}. 
\par
Our investigations in this direction suggest that
the validity of approximations of the firing densities in the
presence of memory effects by the FPT densities of Markov type is
clearly depending on the magnitude of the correlation time.
Hence, an object of our present research is the investigation of those models whose
asymptotic behavior
 becomes increasingly similar as the correlation time 
$\vartheta$ grows larger.

%%%%%%%%%%%%%%%%%%%%%%%%%%%%%%%%%%%%%%%%%%%%%%%%%%%%%%%%%

% FIGURE REFERENCE \ref{f:samplefig}

%\bfg{4}
%\begin{picture}(,)

% ATTENZIONE:
% NON INCLUDERE IL FILE .ps O .eps
% SPEDIRLO COME FILE SEPARATO
% QUI METTERE SOLO LA CAPTION

%\end{picture}
%\fgc{Figure Title}{fig2label}{Figure caption figure caption figure
%caption} \efg

%% ======================================================================= %%
%% FINE FINE FINE FINE FINE FINE FINE FINE FINE FINE FINE FINE FINE FINE F %%
%% ======================================================================= %%

\myline

\end{document}